\documentstyle[11pt,amssymb,amsfonts]{article}
  
\begin{titlepage}

\title{On the index of equivariant Toeplitz operators}

\author{Ulrich Bunke\thanks{Mathematisches Institut, Universit\"at G\"ottingen, Bunsenstr. 3-5, 37073 G\"ottingen, GERMANY, E-mail:bunke@uni-math.gwdg.de}
}

\end{titlepage}

\newcommand{\diag}{{\rm diag}}
\newcommand{\proof}{{\it Proof.$\:\:\:\:$}}

\newcommand{\R}{{\bf R}}
\newcommand{\Mat}{{\rm Mat}}

\newcommand{\Z}{{\bf Z}}

\newcommand{\C}{{\bf C}}

\newcommand{\grad}{{\rm grad}}

\newcommand{\ch}{{\rm\bf ch}}

\newcommand{\cH}{{\cal H}}

\newcommand{\cV}{{\cal V}}

\newcommand{\cW}{{\cal W}}

\newcommand{\sign}{{\rm sign}}

\newcommand{\cN}{{\cal N}}

\def\hB{\hspace*{\fill}$\Box$ \newline\noindent}

\newcommand{\ind}{{\rm index}}

\def\hB{\hspace*{\fill}$\Box$ \\[0.5cm]\noindent}

\newtheorem{prop}{Proposition}[section]
\newtheorem{lem}[prop]{Lemma}
\newtheorem{ddd}[prop]{Definition}
\newtheorem{theorem}[prop]{Theorem}
\newtheorem{kor}[prop]{Corollary}
\newtheorem{ass}[prop]{Assumption}

\begin{document}
\maketitle
\tableofcontents
        \section{Introduction}

The goal of the present paper is to understand aspects of the recent
preprint \cite{nestradulescu99}. While this paper studies traces of commutators
the present paper concentrates on the index theoretic aspects.
This allows for studying the index Toeplitz operators under quite general
assumptions. The basic observation (compare \cite{higson95})
is that the index of the Toeplitz
operator is equal to the index of an associated Callias type operator,
i.e. a Dirac operator with potential.
Callias type operators were thoroughly studied in \cite{bunke9381}.
In particular, their index is very accessible to computation.
In the present note we show how to extend all that to the equivariant case.

\section{The index of Toeplitz operators}

In this Section we review the non-equivariant situation.
The main ideas can be traced back to \cite{higson95}.

Let $(M,g)$ be a complete Riemannian manifold and $E\rightarrow M$ be
a Dirac bundle which may be $\Z_2$-graded or ungraded. The associated
Dirac operator is an unbounded essentially selfadjoint operator on the Hilbert space
$H:=L^2(M,E)$  with domain $C_c^\infty(M,E)$.
\begin{ass}\label{dspec}
We assume that $0$ is an isolated point of the spectrum of $D$.
\end{ass}
Let $\cH$ denote the kernel of $D$, and let $P$ be the orthogonal projection
onto $\cH$. In the graded case we have an orthogonal splitting
$\cH=\cH^+\oplus \cH^-$, $P=P^++P^-$.

Let $C^\infty_g(M)$ denote the space of all bounded smooth functions
such that $df$ vanishes at infinity of $M$. The commutative
$C^*$-algebra $C_g(M)$ is defined as the closure of $C^\infty_g(M)$
inside $C(M)$. Furthermore, let $C_0(M)$ denote the closure of
$C^\infty_c(M)$ in $C(M)$.

If $f\in C(M)$, then $M_f$ denotes the multiplication operator
on $H$ induced by $f$.
\begin{ddd}
We define the following operators on $\cH$, resp. $\cH^{\pm}$:
$$T_f:=PM_fP,\quad T^\pm_f:=P^\pm M_f P^\pm\ .$$
\end{ddd}
We extend this definition to matrix valued functions
$f\in \Mat(n,C(M))$ such that $M_f\in B(H\otimes \C^n)$,
$T_f\in B(\cH\otimes \C^n)$, etc.
\begin{lem}\label{lll1}
If $f\in C_0(M)$, then $T_f$ is compact.
\end{lem}
\proof
If $f\in C_c^\infty(M)$, then $M_f P$ is compact
by Rellich's Lemma since $\cH\subset H^1$, where
the latter space is the Sobolev space defined as the domain
of the closure of the elliptic operator $D$.
The map $f\mapsto M_f$ is continuous. Since the space of
compact operators is norm-closed we conclude that $M_f P$
is compact for all $f\in C_0(M)$. \hB
\begin{lem}\label{lll2}
If $f\in C_g(M)$, then $[P,M_f]$ is compact.
\end{lem}
\proof
It is here where we use the assumption that $0$ is an isolated
point of $\sigma(D)$. Let $B\subset \C$ be a small ball arround
zero such that $\sigma(D)\cap B=\{0\}$. For $\lambda\not\in\sigma(D)$
let $R_D(\lambda):=(\lambda-D)^{-1}$ denote the resolvent.
By holomorphic function calculus we can write
$$P=\frac{1}{2\pi\imath} \int_{\partial B} R_D(\lambda) d\lambda\ .$$
Let $f\in C^\infty_g(M)$. Then we have
\begin{eqnarray*}
[R_D(\lambda),M_f]&=&R_D(\lambda)[D,M_f]R_D(\lambda)\\
&=&R_D(\lambda)c(\grad f) R_D(\lambda)\ ,
\end{eqnarray*}
where $c$ denotes Clifford multiplication.
Since $\grad f$ vanishes at infinity we conclude that
$c(\grad f) R_D(\lambda)$ is compact. In fact, we can approximate
$\grad f$ uniformly by vector fields $X$ of compact support, and
$X R_D(\lambda)$ is again compact by Rellich's Lemma since
$R_D(\lambda)$ maps $H$ continuously to $H^1$.
Thus $[R_D(\lambda),M_f]$ is compact, too.
Again refering to the norm-closedness of the space of compact operators on $H$
we see that
$$[P,M_f]=\frac{1}{2\pi\imath} \int_{\partial B} [R_D(\lambda),M_f]\: d\lambda$$
is compact for all $f\in C^\infty_g(M)$, and hence for all $f\in C_g(M)$. \hB
Note that Lemma \ref{lll1} and \ref{lll2} extend to matrix valued
functions.

The underlying topological space of $M$ can be considered as the spectrum of the
commutative $C^*$-algebra $C_0(M)$. The exact sequence of $C^*$-algebras
$$0\rightarrow C_0(M)\rightarrow C_g(M)\rightarrow C(\partial_g M)\rightarrow 0$$
defines the compactification of $M$ by the Higson corona $\partial_g M$,
where $\partial_g M$ is the spectrum of the commutative $C^*$-algebra
$C(\partial_g M):=C_g(M)/C_0(M)$.

We now consider a continuous unitary matrix valued function
$F:\partial_g M \rightarrow U(n)$, or equivalently,
a unitary $F\in U(n,C(\partial_g M))$. In the ungraded case
we assume in addition that $F$ is selfadjoint.
Let $f,g\in \Mat(n,C_g(M))$ be extensions of
$F$ and $F^{-1}$ to $M$.
\begin{lem}\label{lll3}
The operator $T_{f}$ is a Fredholm operator. Its index
only depends on $F$, where we define the index in the graded
case as $\ind(T_{f}):=\ind(T^+_{f})-\ind(T^-_{f})$.
\end{lem}
\proof
A parametrix of $T_{f}$ is given by $T_{g}$.
In fact, if "$\sim$" denotes equality modulo compact operators, then
we have by Lemma \ref{lll2} that
$T_{f}T_{g}\sim T_{fg}=1+T_{fg-1}$.
Since $fg-1\in \Mat(n,C_0(M))$ we have by
Lemma \ref{lll1} that $T_{fg-1}\sim 0$.
Similarly we show that $T_{g}T_{f}\sim 1$.
If $f,f^\prime$  are two extensions of $F$,
then $f-f^\prime\in \Mat(n,C_0(M))$. By Lemma
\ref{lll1} we conclude that $T_{f}\sim T_{f^\prime}$
and hence equality of the indices.\hB

The goal of the present section is to relate the index
of $T_{f}$ with the index of the Callias-type operator
constructed in \cite{bunke9381} from $D$ and $f,g$.
The Theorems 2.9 and 2.16 of \cite{bunke9381} reduce the
computation of the index of the Callias type operator
and thus of $T_{f}$ to an application of the Atiyah-Singer
Index theorem for elliptic differential operators on
closed manifolds. In the ungraded case  \cite{bunke9381}, 
Prop. 2.8 implies that $\ind(T_{f})=0$, since
we assume that $\sigma(D)$ has a gap.

In the ungraded case we define the Callias-type operator
$C:=D + i M_{f}$ on $E\otimes \C^n$ and put $\epsilon:=1$.
In the graded case we define $\epsilon:=-1$ and
$C:=D + \diag(-M^+_{g},M^-_{f})$ on $E^+\otimes \C^n \oplus E^-\otimes \C^n$.
As shown in Sec. 2 of \cite{bunke9381} these Callias-type
operators have a well-defined index
$$\ind(C):=\dim\ker_{L^2}(C)-\dim\ker_{L^2}(C^*)\ .$$
Our main result is
\begin{prop}\label{po1}
$$\ind(T_{f})=\epsilon \:\ind(C)$$
\end{prop}
\proof
We first consider the graded case.
Let $H^1$ be the domain of the closure of
(the extension of) $D$  (from $E$ to $E\otimes \C^n)$
and put $H:=L^2(M,E\otimes \C^n)$.
Then $C:H^1\mapsto H$ is Fredholm in the usual sense.
Set $\Phi:=\diag(-M_{g},M_{f})$ and $Q:=1-P$.
$P$ and $Q$ act on $H$ as well as on $H^1$.
It follows from Lemma \ref{lll3} that
$Q\Phi P$ and $P\Phi Q$ are compact operators
from $H^1$ to $H$. Therefore we can replace $C$ by
$PCP\oplus QCQ$ without changing the index.
Using the homotopy $QC_tQ$, $C_t:=D + t\Phi$
we can deform $QCQ$ to the invertible operator
$QC_0Q$ through Fredholm operators. Here we use the fact
that $C_t$ is Fredholm for all $t>0$.
We conclude that
\begin{eqnarray*}
\ind(C)&=&\ind(PCP)
=\ind(P\Phi P)\\
&=&\ind(T^+_{g})+\ind(T^-_{f})
=-(\ind(T^+_{f})-\ind(T^-_{f}))\\
&=&-\ind(T_{f})\ .
\end{eqnarray*}
We now come to the ungraded case. We again have
$$\ind(C)=\ind(PCP)+\ind(QCQ)=\ind(PCP)=\ind(T_{f})\ ,$$
since $QCQ$ can be deformed to the invertible operator $QC_0Q$,
where $C_t:=D + i t M_{f}$. \hB

As explained above we conclude with \cite{bunke9381}, Prop. 2.8 that
\begin{kor}\label{ko1}
In the ungraded case we have $\ind(T_{f})=0$.
\end{kor}

\section{$\Gamma$-equivariant Toeplitz operators}

Let $\tilde M\rightarrow M$ be a Galois cover with group of deck
transformations $\Gamma$. We reserve the symbol "$\:\tilde{}\:$" to
denote lifts of various objects to $\Gamma$-coverings.
\begin{ass}
We can choose a cut-off function $\chi^\Gamma\in C_g^\infty(\tilde M)$
such that $\sum_{\gamma\in\Gamma} \gamma^* \chi^\Gamma=1$,
$\sharp\{\gamma\in\Gamma\:|\: \gamma^*\chi^\Gamma \chi^\Gamma\not=0\}<\infty$.
\end{ass}
In our motivating example $\tilde M$ is a symmetric space
of rank one such that $\Gamma$ is a convex-cocompact group of isometries.
In this case $\chi^\Gamma$ exists by \cite{bunkeolbrich982}, Lemma 6.4.
\begin{ass}\label{gdspec}
We assume that $0$ is an isolated point of the spectrum of $\tilde D$.
\end{ass}
If $\tilde D$ is a homogeneous Dirac operator on a symmetric space
$\tilde M$ of rank one and $M$ is a quotient by a convex cocompact subgroup,
then by the result of \cite{bunkeolbrich982} we have $\sigma_{ess}(D)=\sigma_{ess}(\tilde D)$.
In particular, if $\tilde D$ is one of the Dirac operators constructed
by \cite{atiyahschmid77} in order to realize the representations of the
discrete series, then the assumptions \ref{dspec} and \ref{gdspec}
are satisfied.

We consider the Hilbert space $\tilde H:=\L^2(\tilde M,\tilde E)$
which carries an unitary representation of $\Gamma$.
Let $B_\Gamma:={}^\Gamma B(\tilde H)$ denote the $\Gamma$-equivariant
bounded operators on $\tilde H$.
Using a fundamental domain $F\subset \tilde M$ we can write
$\tilde H=L^2(\Gamma)\otimes H$ in a $\Gamma$-equivariant way.
Let $\cN(\Gamma)\subset B(L^2(\Gamma))$
be the group von Neumann algebra of all operators commuting with left
translations. Then $B_\Gamma=\cN(\Gamma)\otimes B(H)$.
Let $K_\Gamma\subset B_\Gamma$ be the ideal of $\Gamma$-compact
operators corresponding to
$\cN(\Gamma)\otimes K(H)$. An operator $A\in B_\Gamma$ is called $\Gamma$-Fredholm if it is invertible modulo $K_\Gamma$.
If $A$ is $\Gamma$-Fredholm, then its index is an element
$\ind(A)\in K_0(K_\Gamma)=K_0(\cN(\Gamma))$. The normalized trace $\tau$
on the $II_1$-factor $\cN(\Gamma)$ induces a homomorphism $\tau: K_0(\cN(\Gamma))\rightarrow \R$. We define the $\Gamma$-index of $A$ by
$\ind_\Gamma(A):=\tau(\ind(A))$.
\begin{lem}\label{lll35}
If $A\in B_\Gamma$ and $M_{\chi^\Gamma} A$ is compact, then
$A$ is $\Gamma$-compact.
\end{lem}
\proof
We make the isomorphism 
$I:B_\Gamma \stackrel{\sim}{\rightarrow} \cN(\Gamma)\otimes B(H)$ explicit. First we identify
$H=L^2(F,\tilde E_{|F})$. Then $i:\tilde H \stackrel{\sim}{\rightarrow} L^2(\Gamma)\otimes L^2(F,\tilde E_{|F})$ is given by $i(\phi):= \sum_{\gamma\in\Gamma} \delta_\gamma \otimes
\chi_F \gamma^{-1} \phi$, where $\delta_\gamma(\gamma^\prime)$
is zero for $\gamma\not=\gamma^\prime$ and $1$ in the remaining case, and
$\chi_F$ denotes the characteristic function of $F$.
The inverse of this identification is given by
$i^{-1}(\sum_{\gamma\in\Gamma}\delta_\gamma \otimes \phi_\gamma)= 
\sum_{\gamma\in\Gamma} \gamma \phi_\gamma$. 
We now compute $I(A)=i\circ A\circ i^{-1}$
\begin{eqnarray*}
I(A)(\sum_{\gamma\in \Gamma}\delta_\gamma \otimes \phi_\gamma)&=&
 \sum_{\gamma,\gamma^\prime \in \Gamma}\delta_{\gamma^\prime} \otimes \chi_F (\gamma^\prime)^{-1} \gamma A \phi_\gamma\\
&=&\sum_{\gamma,\gamma^\prime \in \Gamma}\delta_{\gamma \gamma^\prime} \otimes \chi_F (\gamma^\prime)^{-1}   A \phi_\gamma\\
&=&(\sum_{\gamma^\prime\in\Gamma}
R(\gamma^\prime)\otimes A_{\gamma^\prime})( \sum_{\gamma\in \Gamma}\delta_\gamma \otimes \phi_\gamma)\ ,
\end{eqnarray*}
thus $I(A)=\sum_{\gamma \in\Gamma}
R(\gamma )\otimes A_{\gamma}$, where
$A_\gamma:= \chi_F \gamma^{-1} A \chi_F$ and $R(\gamma)\in\cN(\gamma)$
is the right translation by $\gamma$. 
If $M_{\chi^\Gamma} A$ is compact, then so is $M_{\gamma^*\chi^\Gamma} A$
for all $\gamma\in\Gamma$. For all
$\gamma\in\Gamma$ we see that
$\chi_F \gamma^{-1} A \chi_F = \sum_{\gamma^\prime\in\Gamma}  \chi_F \gamma^{-1} M_{(\gamma^\prime)^* \chi^\Gamma} A \chi_F$
is compact since  $\chi_F \gamma^{-1} M_{(\gamma^\prime)^* \chi^\Gamma}\not=0$
for at most finitely many $\gamma^\prime\in\Gamma$.
Thus $A$ is $\Gamma$-compact. \hB

If $f\in C_g(M)$, then we have the multiplication operator
$M_{\tilde f}\in {}^\Gamma B_\Gamma$. Let $\tilde \cH$
denote the kernel of $\tilde D$ and $\tilde P\in {}^\Gamma B(\tilde H)$
the orthogonal projection onto $\tilde \cH$.
\begin{ddd}
We define $\tilde T_f:=\tilde P M_{\tilde f}\tilde P\in B_\Gamma$.
\end{ddd}
\begin{lem}\label{lll4}
If $f\in C_0(M)$, then $\tilde T_f$ is $\Gamma$-compact.
\end{lem}
\proof
We have
$$M_{\chi^\Gamma} \tilde T_f  = M_{\chi^\Gamma} \tilde P M_{\tilde f} \tilde P
\stackrel{\rm Lemma \ref{lll2}}{\sim}  \tilde P  M_{\chi^\Gamma \tilde f}
\stackrel{\rm Lemma \ref{lll1}}{\sim} 0$$
since $\chi^\Gamma \tilde f\in C_0(\tilde M)$.
The assertion now follows from Lemma \ref{lll35}.
\hB
\begin{lem}\label{lll5}
If $f\in C_g(M)$, then $[\tilde P,M_{\tilde f}]$ is $\Gamma$-compact.
\end{lem}
\proof
We employ the same method as in the proof of Lemma \ref{lll2} using
$$M_{\chi^\Gamma} \tilde R(\lambda) c(\grad\tilde f) \tilde R(\lambda)
 \sim
 \tilde R(\lambda) M_{\chi^\Gamma} c(\grad\tilde f) \tilde R(\lambda)
\sim 0$$
since $M_{\chi^\Gamma} c(\grad\tilde f)$ vanishes at infinity of $\tilde M$.
The assertion now follows from Lemma \ref{lll35}.
\hB

Consider $F\in U(n,C(\partial_g M))$ which is
selfadjoint in the ungraded case. Let $f,g\in
\Mat(n,C_g(M))$ be lifts of $F, F^{-1}$.
\begin{lem}\label{lll6}
The operator $\tilde T_{f}$ is $\Gamma$-Fredholm.
Its index only depends on $F$.
\end{lem}
\proof
The proof is analogous to the proof of Lemma  \ref{lll3}.
One has to replace "compact" by "$\Gamma$-compact
and applies Lemma \ref{lll4} and \ref{lll5} instead of
Lemma \ref{lll1} and \ref{lll2}. \hB  

The following theorem is the main result of the present paper.
\begin{theorem}\label{main}
$$\ind_\Gamma(\tilde T_{f}) = \ind(T_{f})\ . $$
\end{theorem}
\proof
The lifts $\tilde f, \tilde g$ give rise to a $\Gamma$-equivariant
Callias type operator $\tilde C$. In the first step we show that
$\tilde C$ is $\Gamma$-Fredholm, and that its index coincides
with the index of $\tilde T_{f}$ (up to the sign $\epsilon$).
In \cite{bunke9381}, Sec.2 the computation of the index of $C$ was reduced to
the computation of the index of an elliptic differential operator $R$
on a closed manifold $S^1\times N$ using a relative index theorem and a cut-and-past
procedure. Doing this cut-and-paste procedure equivariantly in the second step
we reduce the computation of the index of $\tilde C$ to the computation
of the index of the lift $\tilde R$ of $R$ to a certain cover $S^1\times \tilde N$. In the third and final step we apply the
Atiyah's index theorem for coverings in order
to conclude that the $\Gamma$-index of $\tilde R$ coincides with the index
of $R$.

We form the flat bundle of
von Neumann algebras
$\cV:=\tilde M\times_\Gamma \cN(\Gamma)$
and let $D_\cV$ be the $\cN(\Gamma)$-equivariant twisted Dirac operator
on $E\otimes \cN(\Gamma)$. We now form the $\cN(\Gamma)$-equivariant
Callias type operators
$\hat C:=D_\cV + i M_{f}$ on $E\otimes \cV\otimes \C^n$ in the ungraded and
$\hat C:=D_\cV + \diag(-M^+_{g},M^-_{f})$ on $E^+\otimes \cV\otimes \C^n\oplus
E^-\otimes \cV\otimes \C^n$ in the graded case.
Combining \cite{bunke9381}, Lemma 2.6, 2.14, and the proof of Lemma 1.18
we show that the operator $\hat C$ is invertible at infinity (see
\cite{bunke9381}, Ass. 1).
It follows that $\hat C$ induces a Fredholm operator
between the Hilbert-$\cN(\Gamma)$ modules $H^1(M,E\otimes \cV\otimes \C^n)$ and
$L^2(M,E\otimes \cV\otimes \C^n)$. The tensor products over $\cN(\Gamma)$ of these modules
with $L^2(\Gamma)$ identify with
$H^1(\tilde M,\tilde E\otimes \C^n)$ and
$L^2(\tilde M,\tilde E\otimes \C^n)$, respectively.
The operator $\hat C$ gives rise to
the $\Gamma$-equivariant Callias type operator
$\tilde C$ which is just the lift of $C$. In particular,  we see that
$\tilde C$ is $\Gamma$-Fredholm and $\ind(\hat C)=\ind (\tilde C)$.
We can now apply exactly the same argument as in the proof of Proposition
\ref{po1} in order to show that
$\ind(\tilde C)=\epsilon \:\ind(\tilde T_{f})$,
replacing compactness and Fredholm by the corresponding $\Gamma$-equivariant notions.
This ends the first step of the proof.

We now come to the second step. In the ungraded case we can repeat
the argument of the proof of \cite{bunke9381}, Prop. 2.8  in order to see that
$\ind(\hat C)=0$ since by assumption there is a gap in the spectrum
of $D_\cV$ (note that the spectrum of $D_\cV$ coincides with that of $\tilde D$).
Thus $0=\ind_\Gamma(\tilde T_{f})$. Since $\ind(T_{f})=0$ by Corollary
\ref{ko1} we obtain the assertion on the theorem in the ungraded case.
It remains to consider the graded case.
Doing the construction \cite{bunke9381}, 2.4.2 with
$D$ and $D_\cV$ at the same time we arrive at a Dirac operator
$R$ and its twist $R_\cW$ over a compact manifold $S^1\times N$
such that $\ind(C)=\ind(R)$ and $\ind(\hat C)=\ind(R_\cW)$.
Here $N$ is a certain closed hypersurface of $M$,
$R$ is associated to a Dirac bundle $L\mapsto S^1\times N$,
and $\cW$ is the pull-back to $S^1\times N$ of the restriction of $\cV$ to $N$.
Let $\tilde N$ be the restriction of the cover $\tilde M\rightarrow M$
to $N$. The tensor products over $\cN(\Gamma)$ of $H^1(S^1\times N,L\otimes \cW)$, $L^2(S^1\times N,L\otimes \cW)$
with $L^2(\Gamma)$ identify with $H^1(S^1\times \tilde N,\tilde L)$ and
$L^2(S^1\times \tilde N,\tilde L)$, respectively. The operator $R_\cW$ induces
the $\Gamma$-equivariant Dirac operator $\tilde R$ on $\tilde L$ which is just the lift
of $R$. We have $\ind(\tilde R)=\ind(R_\cW)$. This accomplishes the second step.

In the last step we apply Atiyah's index theorem for coverings
\cite{atiyah76} in order to conclude that
$\ind_\Gamma(\tilde R)=\ind(R)$. This implies
$\ind_\Gamma(\tilde T_{f})=\ind(T_{f})$ by
Proposition \ref{po1} and the first two steps.
\hB

\section{Examples}

We first consider the two-dimensional example.
Let $\tilde M$ be the hyperbolic plane and $\Gamma$ be a convex
cocompact subgroup of the group of isometries of $\tilde M$.
The geodesic boundary $\partial \tilde M$ can be decomposed into
a limit set $\Lambda$ and its complement $\Omega$.
The group $\Gamma$ acts freely and properly on $\tilde M\cup \Omega$,
and the compact manifold with boundary
$\bar M:=\Gamma\backslash \tilde M\cup \Omega$ is the geodesic
compactification of $M:=\Gamma\backslash \tilde M$.
The boundary $B:=\partial \bar M$ is a finite union of circles
$\Gamma\backslash \Omega$.
There is a natural projection of the Higson corona $\partial_g M$ to
$B$. Thus any $U(1)$-valued function $F$ on $B$ can be lifted to $\partial_g M$,
and we will denote this lift by the same symbol.

Note that $M$ is a complex manifold. Let $K$ be the canonical bundle
of $M$. We fix $k\in\Z$ and consider the graded Dirac operator $D=\bar \partial +(\bar \partial)^*$ on
$E=K^k\oplus K^{k-1}$, where $\bar \partial : C^\infty(M,K^k)\rightarrow C^\infty(M,K^{k-1})$
is the Dolbeault operator.

The complex structure fixes an orientation of $M$
which induces an orientation of $B$.
We now use the notation of \cite{bunke9381}, 2.16. The Dirac operator
$D_N$ is just $i\frac{\partial}{\partial t}$ on any component
of $B$, where $t$ is the coordinate of $S^1$ compatible with the orientation.
The index $\ind(T_{f})$ is minus the  spectral
flow of the family connecting $D_N$ and $F^*D_NF$,
and this is equal to the total winding number $n(F)$ of $F$,
i.e.  $$\ind(T_{f})=-\frac{1}{2\pi\imath}\int_B F^{-1}dF\ .$$

Let $P_k$ denote the projection onto the space of holomorphic square integrable sections of $\tilde K^k$. Note that $P_0=0$. For $k>0$ ($k<0$) the range of
$P_k$ is the holomorphic (antiholomorphic) discrete series representation
of $PSL(2,R)$, the orientation-preserving isometry group of $\tilde M$.
Let $\tilde T^k_f:=P_k\tilde fP_k$ be the Toeplitz compressions.
Then by the computation above and Theorem \ref{main} for $k\not=0$ we have
$\ind_\Gamma(\tilde T^k_f)=-\sign(k) n(F)$.

This has the following higher-dimensional generalization.
Let $\tilde M:=Spin(1,2n)/Spin(2n)$
be the real hyperbolic space of dimension $2n$
and $\Gamma\subset Spin(1,2n)$ be a convex cocompact subgroup.
Let $\tilde S$ be the spinor bundle of $\tilde M$
and $\tilde V$ be any further $Spin(1,2n)$-homogeneous bundle.
We put $\tilde E:=\tilde S\otimes \tilde V$ and
let $\tilde D$ be the associated Dirac operator.
We assume that $\tilde V$ is such that $0$ is an
isolated point of the spectrum of $\tilde D$.
In this case $\ker \tilde D$ decomposes into a finite sum of
discrete series representations of $Spin(1,2n)$.

We again have a decomposition of the geodesic boundary of $\tilde M$
into a limit set and a domain of discontinuity $\Omega$.
The locally symmetric space $M:=\Gamma\backslash \tilde M$ can
be compactified by adjoining the boundary $B:=\Gamma\backslash \Omega$.
The topology of $B$ can be quite complicated. Since $Spin(1,2n)$
acts on the sphere $\partial \tilde M$ by orientation-preserving
conformal transformations $B$ admits a locally conformally
flat structure. In particular,
all Pontrjagin classes of $TB$ and all associated bundles vanish.

Again we have a natural map from the Higson corona
of $M$ to $B$. Let $F:B\rightarrow U(m)$ be a continuous function
and $f,g\in Mat(m,C_g(M))$ extensions of $F,F^{-1}$.
The function $F$ represents an element $[F]$ in $K^1(B)$.
Let $\ch([F])\in H^{odd}(B,\R)$ be the Chern character
of $[F]$. We define the degree of $F$ by
$\deg(F):=\ch_{2n-1}([F])([B])$. Here $[B]$ is the orientation of $B$
as the boundary of $M$, where the orientation of $M$ is determined by
the $\Z_2$-grading of $S$.
\begin{prop}
$\ind_\Gamma(\tilde T_f)=\ind(T_f)=-\dim(V)\deg(F)$.
\end{prop}
\proof
We have to compute the index of the Callias type
operator $C$ on $M$ which is associated to $f,g$.
By \cite{bunke9381}, Thm. 2.16 it is equal to the index
of the Dirac operator $D_L$ on $S^1\times B$
twisted with a bundle $L=L_F\otimes L_V$.
Here $L_V$ is the restriction of $V$ to $B$. The bundle
$L_F$ is obtained from the trivial bundle $[0,1]\times B \times \C^m$
by glueing $(1,b,v)$ with $(0,b,F(b)v)$, $b\in B$, $v\in \C^m$.
The bundle $L_V$ is associated to the tangent bundle of $B$
and thus has vanishing Chern classes. Further note
that ${\bf \hat{A}}(T(S^1\times N))=1$. The index theorem for
twisted Dirac operators thus gives
$$\ind(D_L)=\dim(V)\ch(L_F)_{[2n]}([S^1\times N])=
\dim(V)\ch_{2n-1}([F])([B])\ .$$
This finishes the proof of the proposition.\hB

In the situation above we know that $\ker D$ is infinite-dimensional
by \cite{bunkeolbrich982}. The proposition above would give
an alternative index-theoretic proof of this fact.

\bibliographystyle{plain}

\begin{thebibliography}{1}

\bibitem{atiyahschmid77}
M.~Atiyah and W.~Schmid.
\newblock A geometric construction of the discrete series for semisimple Lie
  groups.
\newblock {\em Invent. Math.}, 42(1977), 1--62.

\bibitem{atiyah76}
M.~F. Atiyah.
\newblock Elliptic operators, discrete groups and von Neumann algebras.
\newblock {\em Asterisque}, 32(1976), 43--72.

\bibitem{bunke9381}
U.~Bunke.
\newblock A {\bf K}-theoretic relative index theorem and Callias-type operators.
\newblock {\em Math. Ann.}, 303(1995), 241--279.

\bibitem{bunkeolbrich982}
U.~Bunke and M.~Olbrich.
\newblock The spectrum of Kleinian manifolds.
\newblock To appear in  J. Funct. Anal., Preprint available at
  http://www.uni-math.gwdg.de/bunke/spzerl.dvi.

\bibitem{higson95}
E.~Guenter and N.~Higson.
\newblock A note on Toeplitz operators.
\newblock {\em Int. J. Math.}, 7(1996), 501--513.


\bibitem{nestradulescu99}
R.~Nest and F.~Radulescu.
\newblock Index of $\Gamma$-equivariant Toeplitz operators.
\newblock Preprint 1999 : math.OA/9911042.


\end{thebibliography}

\end{document}